\documentclass [12pt,a4paper]{article}
\usepackage[T1]{fontenc}
\usepackage[ansinew]{inputenc}
\usepackage{amsthm,amsmath,amssymb}
\usepackage{graphicx}
\usepackage{cite}

\newcommand{\R}{\mathbb R}
\newcommand{\s}{\mathbb S}

\newcommand{\eps}{\varepsilon}

\newtheorem{thm}{Theorem}[section]
\newtheorem{cor}{Corollary}[section]
\newtheorem{lem}{Lemma}[section]
\newtheorem{prop}{Proposition}[section]

\newtheorem{rem}{Remark}[section]

\begin{document}

\title{Global existence for the generalized two-component   Hunter-Saxton system}

\author{{\sc Hao Wu} \footnote{Shanghai Key Laboratory for Contemporary Applied Mathematics and School of Mathematical Sciences, Fudan
University, 220 Han Dan Road, 200433 Shanghai, China. Email:
\emph{haowufd@yahoo.com}} \ and {\sc Marcus Wunsch}
\footnote{Research Institute for Mathematical Sciences,
 Kyoto University, Kyoto 606-8502, Japan, Email: \emph{mwunsch@kurims.kyoto-u.ac.jp}} }
\date{August 7, 2010}

 \maketitle

\begin{abstract}
 \noindent
We study the global existence of solutions to a two-component
generalized  Hunter-Saxton system in the periodic setting. We first
prove a persistence result of the solutions. Then for some
particular choices of the parameters $(\alpha, \kappa)$, we show the
precise blow-up scenarios and the existence of global solutions to
the generalized Hunter-Saxton system under proper assumptions on the
initial data. This significantly improves recent results obtained in
\cite{Wun09,Wun10}.
\end{abstract}

{\it Keywords:} Generalized Hunter-Saxton system; global
existence, blow-up scenario

{\it 2010 Mathematics Subject Classification}
35B10, 
35B65; 
35Q35; 
35Q85 
\noindent


\section{Introduction}
\setcounter{equation}{0}

\noindent
In this paper, we shall investigate the global existence of solutions to
 the periodic boundary value problem for a two-component family of evolutionary
 systems modeling fluid convection and stretching in one space dimension, \\[0.1cm]
\begin{equation} \label{gen:hss}
\begin{cases}
\partial_t \; m(t,x) + \underbrace{u \partial_x \;m}_{{\rm convection}} + \;\; (1 - \alpha)
\underbrace{\partial_x u \;m}_{{\rm stretching}} +
\underbrace{\kappa \; \rho \;\partial_x \rho}_{{\rm coupling}} = 0,
\\[0.1cm]
 m(t,x) = -\partial_{xx}^2 \; u(t,x), \\[0.1cm]
\partial_t\;\rho + \underbrace{u\partial_x \; \rho}_{ {\rm convection}} = \alpha  \; u_x \rho, \\[0.5cm]
m(0,x) = m^0(x), \; \rho(0,x) = \rho^0(x), \quad x \in \s \simeq \R/\mathbb Z,\\[0.1cm]
\end{cases}
\end{equation}
where $(1 - \alpha) \in \R$ is the ratio of stretching to
convection, and $\kappa$ denotes a real dimensionless constant that
 measures the impact of the coupling.
\par
This system was first studied in full generality by {\sc Wunsch}
\cite{Wun10}, where it was coined the {\it generalized Hunter-Saxton
system}, since for $(\alpha,\kappa) = (-1,\pm 1)$ it becomes the
Hunter-Saxton system \cite{Wun09}. The latter is a particular case
of the Gurevich-Zybin system pertaining to nonlinear one-dimensional
dynamics of dark matter as well as nonlinear ion-acoustic waves (cf.
{\sc Pavlov} \cite{pavlov} and the references therein).
\par
It was noted by {\sc Constantin \& Ivanov} \cite{ci} that the
Hunter-Saxton system allows for peakon solutions; moreover, {\sc
Lenells  \& Lechtenfeld} \cite{ll} showed that it can be interpreted as
the Euler equation on the superconformal algebra of contact vector
fields on the $1 \big| 2$-dimensional supercircle, which is in
accordance with the by now well-known geometric interpretation of
the Hunter-Saxton  equation as the geodesic flow of the
right-invariant $\dot H^1(\s)$ metric on the space of
orientation-preserving circle diffeomorphisms modulo rigid rotations
\cite{km,ll, len08,len:weak} (see also \cite{ck:geod,
ck:integ,ck,kour,misiolek} for related geodesic flow equations).

The two-component
Hunter-Saxton system is a generalization of the Hunter-Saxton
equation modeling the propagation of weakly nonlinear orientation
waves in a massive nematic liquid crystal (see {\sc Hunter \&
Saxton} \cite{h-s} for a derivation, and also
\cite{bss,bc,bhr,tiglay,yin}), since the former obviously reduces to
the latter if the initial datum $\rho^0$ is chosen to vanish
identically. It turns out that if this choice is made for arbitrary
$\alpha \in \R$, one arrives at the generalized Proudman-Johnson
equation \cite{pj,okamoto,oo,oz,cw,w} with parameter $a = \alpha -
1$. Through this link, the family of systems \eqref{gen:hss} also
bridges the rich theories for the Burgers equation \cite{burgers}
($\alpha = -2$) and the axisymmetric Euler flow in $\R^d$
\cite{okamoto, st} if $\alpha = 2/(d - 1)$. We also remark that if
one sets $\rho = \sqrt{-1} \; u_x$ and $\kappa = - \alpha$, the
system \eqref{gen:hss} decouples to give, once again, the
generalized Proudman-Johnson equation \cite{okamoto, chow, w} with
parameter $a = 2\alpha - 1$.
Other important special cases of the generalized Hunter-Saxton
system \eqref{gen:hss} include the inviscid K\'arm\'an-Batchelor
flow \cite{chae,co,hl} for $\alpha = -\kappa = 1$, which admits
global strong solutions, and the celebrated Constantin-Lax-Majda
equation \cite{clm} with $\alpha = -\kappa = \infty$, a
one-dimensional model for three-dimensional vorticity dynamics,
which has an abundance of solutions blow-up in finite time.
\par
The Hunter-Saxton system \eqref{gen:hss} with parameters
$(\alpha,\kappa) = (-1,\pm 1)$ is the \textit{short-wave limit},
obtained via the space-time scaling $(x,t) \mapsto (\eps x, \eps t)$
and letting $\eps$ tend to zero in the resulting equation, of the
two-component integrable Camassa-Holm system \cite{ci,ELY07}. This
system, reading as \eqref{gen:hss} with $m$ replaced by $(1 -
\partial_{xx}^2) u$, has recently been the object of intensive study
(see \cite{ci,ELY07,gy,guo,gz,ll,mohajer,mustafa,ZL09}). {\sc
Constantin \& Ivanov} \cite{ci} derived the Camassa-Holm system from
the Green-Naghdi equations, which themselves originate in the
governing equations for water waves \cite{johnson}.
\par
One major motivation for studying systems such as the Camassa-Holm
system or the system \eqref{gen:hss} lies in their potential
exhibition of nonlinear phenomena such as wave-breaking and peaked
traveling waves, which are not inherent to small-amplitude models
but known to exist in the case of the governing equations for water
waves (prior to performing asymptotic expansions in special regimes
like the shallow water regime), cf.  \cite{toland, constantin,
whitham, johnson}. In this context, it is of interest to point out
that peaked solitons are absent among the solitary wave solutions to
the Camassa-Holm system (cf. \cite{mustafa}), while they exist for the Hunter-Saxton
system, see \cite{ci}.
\par
Another reason - and, indeed, the very incentive in \cite{Wun10} and
here - for analyzing the family of systems \eqref{gen:hss}, has its
origin in a paradigm of {\sc Okamoto \& Ohkitani}\cite{oo} that the
convection term can play a positive role in the global existence
problem for hydrodynamically relevant evolution equations (see also \cite{hl,osw}). The quadratic terms in the first component
of \eqref{gen:hss} represent the competition in fluid convection
between nonlinear steepening and amplification due to $(1 -
\alpha)$-dimensional stretching and $\kappa$-dimensional coupling
(cf. \cite{hs}). The stretching parameter $\alpha$ illustrates the inherent importance of the convection term in delaying or depleting finite-time blow-up, while the coupling constant $\kappa$ measures the strength of
the coupling, and has a strong influence on singularity formation or global existence of the solutions.
Recently, {\sc Wunsch} \cite{Wun10}
proved that the first solution component breaks down in finite time
if $(\alpha,\kappa) \in \{-1\}\times \R_-$ and if the initial slope
is large enough; moreover, he demonstrated, for $(\alpha,\kappa) \in
[-1/2,0)\times \R_-$, that a sufficiently negative slope at an
inflection point of $u^0$ will become vertical spontaneously. By
analogy with  the Constantin-Lax-Majda vorticity model equation
\cite{clm}, the case of $\infty$-dimensional stretching and coupling (i.e., $\alpha=\kappa=\infty$)
was shown to lead to catastrophic steepening of the first solution
component $u$ as well.\medskip

\textbf{Outline of results.} The main purpose of this paper is to broaden our
understanding of solutions to \eqref{gen:hss} by proving rigorously
that solutions for some particular cases (e.g., $(\alpha,\kappa) \in
\{-1,0\} \times \R_+$) can be global. In the preliminary Section
\ref{suf:con}, we first recall the local-in-time
well-posedness result of system \eqref{gen:hss} for $(\alpha, \kappa)\in
\mathbb{R}\times \mathbb{R}$ and provide a partial result on the rate of break-down at the origin for $(\alpha, \kappa)\in \{-1\}\times \mathbb{R}_{-}$. In Section 3, we first prove a persistence result of solutions in
$H^s\times H^{s-1}$, $s\geq 2$, for $(\alpha, \kappa)\in
\mathbb{R}\times \mathbb{R}$. In Section 4, we derive
some precise blow-up scenarios for the solutions in the case
$(\alpha,\kappa) \in \{-1\} \times \R_+$. In Section 5, we show the
existence global solutions in $H^s\times H^{s-1}$ (only beyond a
small regularity gap $s\in (2, \frac52])$ under proper assumptions
on the initial data, which replaces the
artificial assumption made in
Section 3 that the gradient of the second solution component $\rho$
be bounded. In Section 6, the global existence of sufficiently regular
solutions when $(\alpha,\kappa) \in \{ 0 \} \times \R_+$ is
discussed. \medskip

\textbf{Notations.} Throughout the paper, $\s = \R / \mathbb Z$
shall denote the unit circle. By $H^r(\s)$, $r\ge 0$, we will
represent the Sobolev spaces of equivalence classes of functions
defined on the unit circle $\s$ which have square-integrable
distributional derivatives up to order $r$. The $H^r(\s)$-norm will
be designated by $\|.\|_{H^r}$ and the norm of a vector $b \in
H^r(\s) \times H^{r-1}(\s)$ will be written as $\|b\|_{H^r\times H^{r-1}}$.
Also, the Lebesgue spaces of order $p \in [1,\infty]$
will be denoted by $L^p(\s)$, and the norms of their elements by $\|f\|_{L^p}$.
Finally,
if $p = 2$, we agree on the convention $\|.\|_{L^2} =: \| .\|$;
moreover, $\langle.,.\rangle := \langle.,.\rangle_{L^2}$ will denote
the $L^2$ inner product. The relation symbol $\lesssim$ stands for
$\leq C$, where $C$ denotes a generic constant.

 \section{Preliminaries} \label{suf:con}
 \setcounter{equation}{0}
  We rewrite the system \eqref{gen:hss} and consider the following
  problem with periodic boundary conditions in the remaining part of the paper:
\begin{eqnarray} \label{hssys}
\begin{cases}
u_{txx} + (1-\alpha)\;u_x\;u_{xx} + u\;u_{xxx} - \kappa \rho\;\rho_x = 0, \;\ t>0,\; x \in \R, \\
\rho_t + u\;\rho_{x} = \alpha u_x \rho, \quad t>0,\; x \in \R, \\
 u(t,x + 1) = u(t,x),\ \ \rho(t,x + 1) = \rho(t,x)\quad t\geq 0,\;x\in \R, \\
 u(0,x) =  u^0(x), \ \ \rho(0,x) =  \rho^0(x),\quad x\in \R.
\end{cases}
\end{eqnarray}
 \begin{rem}
 Integration in space of the $u$-equation in \eqref{hssys} yields
 \begin{equation}
 u_{tx}+uu_{xx}-\frac{\alpha}{2}u_x^2-\frac{\kappa}{2}\rho^2=a(t),\label{HS1}
 \end{equation}
 where the time-dependent integration constant $a(t)$ is determined
 by the periodicity of $u$ to be
 \begin{equation}
 a(t)=-\frac{\kappa}{2}\int_\mathbb{S}\rho^2 dx -\frac{\alpha+2}{2}\int_\mathbb{S}
 u_x^2dx.\label{at}
 \end{equation}
Integrating in space once more, one gets
 \begin{equation}
 u_t+uu_x=\partial_x^{-1}\left(\frac{\kappa}{2}\rho^2+\frac{\alpha+2}{2}u_x^2+a\right)+h(t),\label{HS2}
 \end{equation}
 where $\partial_x^{-1}f(x):=\int_0^x f(y)dy$ and
 $h(t):[0,+\infty)\to\mathbb{R}$ is an arbitrary continuous
 function.
 \end{rem}

 We first recall a local well-posedness result of system \eqref{hssys} (cf. \cite[Theorem 2.1]{Wun10}, see also \cite[Theorem
 4.1]{Wun09} for the special case $\alpha=-1, \kappa=1$).
 \begin{thm}\label{loc}%
 Denote $z=(u, \rho)^{tr}$. Given any $z^0=(u^0, \rho^0)^{tr}\in H^s(\mathbb{S})\times H^{s-1}(\mathbb{S})$, $s\geq
 2$, for $(\alpha, \kappa)\in  \mathbb{R}\times\mathbb{R}$, there exists a maximal life span
 $T=T(\|z^0\|_{H^s\times H^{s-1}})>0$ and a unique solution $z$ to
 system \eqref{hssys} that
 $$z\in C([0, T); H^s(\mathbb{S})\times H^{s-1}(\mathbb{S}))\cap C^1([0,T); H^{s-1}(\mathbb{S})\times
 H^{s-2}(\mathbb{S})).$$
 \end{thm}

 \begin{rem}
 Following the arguments in \cite{yin}, it is possible to show that the maximal existence time $T$ of the solution in Theorem \ref{loc}
  can be chosen independently of the Sobolev order $s$.
 \end{rem}

 \begin{lem}[cf. \cite{Wun10}] \label{conva}
 For $(\alpha, \kappa)\in  \mathbb{R}\times\mathbb{R}$, let
 $(u,\rho)$ be a smooth solution to system \eqref{hssys}.
 Then
 \begin{equation}
 \frac{d}{dt}\; a(t)=-\frac32 \kappa(\alpha+1)\int_\mathbb{S}u_x\rho^2
 dx-\frac{(\alpha+1)(\alpha+2)}{2}\int_\mathbb{S} u_x^3 dx.
 \end{equation}
 In particular, if $(\alpha, \kappa)=\{-1\}\times \mathbb{R}$, then
 $\frac{d}{dt} a(t)=0$, which implies that the system enjoys a conservation law, namely,
 \begin{equation}
  a(t)\equiv a(0)= -\frac{1}{2}\|u^0_x\|^2-\frac{\kappa}{2}\|\rho^0\|^2\label{conaa}
 \end{equation}
 is constant for all $t\geq 0$.
 \end{lem}

In contrast with the cases $(\alpha,\kappa) \in \{ -1,0 \} \times \R_+$ where we shall get global existence  (see the subsequent sections), a slope of singularity of system \eqref{hssys} has been obtained for the case
$(\alpha,\kappa) \in \{ -1 \} \times \R_-$ (cf. \cite[Proposition 3.2]{Wun10}), however, no estimate on the blow-up rate was given there. In what follows,
we provide a partial blow-up result at the origin $x=0$ with blow-up rate for solutions to \eqref{hssys}.

\begin{prop}\label{BL1}
Let $z(t,x)$ be a solution to \eqref{hssys} with parameters $(\alpha,\kappa) \in \{ -1 \} \times \R_-$ and initial datum $z^0 \in H^s \times H^{s - 1}$, $s \ge 2$. In addition, we assume that $u^0$ is odd with $u_x(0)< 0$ and $\rho^0$ is even with $\rho^0(0)=0$,
and that,
moreover,
\begin{equation}
\|u^0_x\|^2+\kappa\|\rho^0\|^2\geq 0.\label{blowupi}
\end{equation}
Then at the origin $x=0$, $u_x(t,0)$ blows up in finite time $T_0$ (time of break-down at the origin).
The blow-up rate of
$u_x(t,0)$ is
$$\lim_{t \rightarrow T_0} \left\{ (T_0 - t) \; u_x(t,0) \right\} = - 2.$$
\end{prop}

\begin{proof}
First, we notice that $u(t,.)$, $\rho(t,.)$ remain odd or even, respectively, due to the algebraic structure of the equations in \eqref{hssys}.
Observe next that
\begin{equation} \label{rh0}
\rho(t,0) = 0
\end{equation}
for all times of existence.
Indeed, one has
$$\frac{\partial}{\partial t} \ \rho (t,0) = - u\rho_x(t,0)  - u_x \rho(t,0).$$
Note that the first term on the right-hand side vanishes since both $u$ and $\rho_x$ are odd.
Together with the assumption on initial data that $\rho^0(0) = 0$, it proves the claim \eqref{rh0}.
\\
Let us now set $$\zeta(t) := u_x(t,0).$$
The resulting ordinary differential equation for the evolution of $\zeta$ reads as follows:
\begin{equation}
\begin{cases}
\frac{d}{dt}\zeta(t) = -\frac{1}{2} \zeta(t)^2 + \frac{\kappa}{2} \ \rho(t,0)^2 +a(t) \stackrel{\eqref{rh0}}{=} -\frac{1}{2} \zeta(t)^2+a(t), \\
\zeta(0) = \zeta^0 = u_x^0(0).
\end{cases}\label{ode}
\end{equation}
When $\alpha=-1$, due to Lemma \ref{conva}, $a(t)$ is a constant. Then by \eqref{blowupi} we know that $a\leq 0$.
Thus, we have
$$ \frac{d}{dt}\; \zeta(t) \leq -\frac{1}{2} \zeta(t)^2,$$
which implies
\begin{equation*}
\zeta(t) \le \frac{2\zeta^0}{2 + \zeta^0\; t}
\end{equation*}
As a consequence,
\begin{equation}
\lim_{t\to T_0}\zeta(t)= -\infty. \label{blp}
\end{equation}
 For $a(0)=-\frac{\kappa}{2}\|\rho^0\|^2 -\frac{1}{2}\|u^0_x\|^2\leq 0$, the following chain of inequality holds
\begin{equation} \label{chain}
 -|a(0)|-1 < \frac{d}{dt}\; \zeta(t)+\frac{1}{2} \zeta(t)^2<|a(0)|+1 .
\end{equation}
Because of \eqref{blp}, there exists a number $\eps \in (0,\frac12)$ and a time $t_\eps$ such that $$\zeta(t)^2 \ge \frac{|a(0)|+1}{\eps}>0,\quad \forall t \in (t_\eps,T_0).$$
Hence
$$-\frac{1}{2} - \eps < \frac{1}{\zeta(t)^2}\;\frac{d}{dt} \zeta(t)  < -\frac{1}{2}+\eps, \quad \forall t \in (t_\eps,T_0),$$
from which we glean, upon integrating from $t > t_\eps$ to $T_0$, that
$$-\frac{1}{2} - \eps < \frac{1}{(T_0 - t)\ \zeta(t)} < -\frac{1}{2}+\eps. $$
We may thus conclude the assertion of the proposition, since $\eps$ was chosen arbitrarily.
\end{proof}

\begin{rem}
It remains an open problem to determine the {\rm first time of break-down},
since the ODE describing the evolution of $m(t) := \inf_x u(t,x)$ is more involved than \eqref{ode}, and double-sided estimates of $\frac{d}{dt}m(t)$ -- as in \eqref{chain} -- would require uniform bounds on $\| \rho(t,.) \|_{L^\infty}$ (cf. \cite{Wun10}). We observe that the rate of break-down we obtained is in accordance the one computed for the Camassa-Holm system \cite{guo}.
\end{rem}

\begin{rem}
For some special cases, the exact blow-up time $T_0$ can be computed. For instance,
\begin{itemize}
\item[(i)] $\| u^0 \|^2 = -\kappa \| \rho^0 \|^2 $, then $a(0) = 0$.
\item[(ii)] $\|u^0_x\|^2 = -\kappa \| \rho^0 \|^2 + 1$, then $a(0) = -\frac12$.
\end{itemize}
In case $(i)$, the explicit solution to \eqref{ode} with $a(0)=0$ reads
$$\zeta (t) = \frac{2\zeta^0}{2 + \zeta^0\; t} < 0.$$
Then the blow-up time is given by $$T_0 = -2/\zeta^0>0.$$
In case $(ii)$, the explicit solution to \eqref{ode} with $a(0)=-\frac12$ reads
\begin{equation} \label{zeta}
\zeta(t) = \tan\left( \arctan (\zeta^0) - \frac{t}{2} \right) < 0.
\end{equation}
Thus the blow-up time of $\zeta(t)=u_x(t,0)$
can be given exactly as $$T_0 = \pi + 2 \arctan(\zeta^0) \; \in\;  (0,\pi).$$
\end{rem}

Even
if the condition \eqref{blowupi} does not hold, we can still construct some solutions that break down at the origin.
\begin{cor}
Let $z(t,x)$ be a solution to \eqref{hssys} with parameters
$(\alpha,\kappa) \in \{ -1 \} \times \R_-$ and initial datum $z^0 \in H^s \times H^{s - 1}$, $s \ge 2$. In addition, we assume that $u^0$ is odd with $u_x(0)< 0$ and $\rho^0$ is even with $\rho^0(0)=0$. Moreover, if, instead of \eqref{blowupi}, we assume that
\begin{equation}
u_x^0(0)<- \sqrt{2\left|-\frac12\|u^0_x\|^2-\frac{\kappa}{2}\|\rho^0\|^2\right|},\label
{blowupj}
\end{equation}
then 
$u_x(t,x)$ blows up at the origin $x = 0$ in finite time.
The blow-up rate of $u$ is
$$\lim_{t \rightarrow T_0} \left\{ (T_0 - t) \; u_x(t,0) \right\} = - 2.$$
\end{cor}
\begin{proof}
Similarly to the proof of Proposition \ref{BL1}, we have
\begin{equation}
\frac{d}{dt}\; \zeta(t) = -\frac{1}{2} \zeta(t)^2+a(0) \leq -\frac{1}{2} \zeta(t)^2+|a(0)|.\label{bj}
\end{equation}
Thus, if \eqref{blowupj} holds, namely, $\zeta(0)<-(2|a(0)|)^\frac12$, then $\zeta(t)<-(2|a(0)|)^\frac12$ for all $t\in [0,T_0)$, where $T_0>0$ is the existence time (ensured by Theorem \ref{loc}). By solving the standard Riccati type inequality \eqref{bj}, it follows that (cf. e.g., \cite{guo})
\begin{equation}
\lim_{t\to T_0} \zeta(t)=-\infty, \quad \text{with} \ \quad 0<T_0< (2|a(0)|)^{-\frac12}\ln\frac{\zeta^0-(2|a(0)|)^\frac12}{\zeta^0+(2|a(0)|)^\frac12}.
\end{equation}
The
computation
of the blow-up rate is now exactly the same as for Proposition \ref{BL1}.
\end{proof}


 \section{Persistence of solutions for $(\alpha, \kappa) \in \mathbb{R}\times\mathbb{R}$}
 \label{Per}
 \setcounter{equation}{0}

In this section, we  consider the question of looking for a suitable bound on
the solutions to \eqref{hssys}, which will ensure that the local
solutions can be extended beyond the maximal existence time given by
the local well-posedness theory in Theorem \ref{loc}.

We first introduce some lemmata that are useful in the subsequent estimates:
\begin{lem} \label{KP}[Kato-Ponce commutator estimate]  Denote  $\Lambda = (1 -
\partial_x^2)^{1/2}$. For $s > 0, \; p\in(1,\infty)$,
\begin{equation}
\| [\Lambda^s,f]\; v\|_{L^p} \lesssim \| f_x \|_{L^\infty} \|
\Lambda^{s-1} v \|_{L^p} +  \| \Lambda^s f\|_{L^p} \| v
\|_{L^\infty}.\nonumber
\end{equation}
\end{lem}
\begin{lem}\label{AL}
If $s>0$, then $H^s\cap L^\infty$ is an algebra. Moreover,
\begin{equation}
\|fg\|_{H^s}\lesssim
\|f\|_{L^\infty}\|g\|_{H^s}+\|f\|_{H^s}\|g\|_{L^\infty}.\nonumber
\end{equation}
\end{lem}

The main result of this section is as follows
\begin{thm} \label{T1}
Suppose  $(\alpha,\kappa)\in \mathbb{R}\times \mathbb{R}$. For any
$z^0=(u^0, \rho^0)^{tr}\in H^s(\mathbb{S})\times
H^{s-1}(\mathbb{S})$, $s\geq
 2$, let $T$ be the existence time of the solution $z=(u,
 \rho)^{tr}$ to system \eqref{hssys} corresponding to $z^0$. If
 there exists a constant $M>0$ such that
 \begin{equation}
 \|u_x(t,\cdot)\|_{L^\infty}+\|\rho(t,
 \cdot)\|_{L^\infty}+\|\rho_x(t,\cdot)\|_{L^\infty}\leq M, \quad \forall t\in [0,
 T),\label{Asm}
 \end{equation}
 then $\|z(t, \cdot)\|_{H^s\times H^{s-1}}$ is bounded on $[0,T)$.
\end{thm}
\begin{proof}
In the proof, we perform only formal calculations which can, however, be justified rigorously
using Friedrichs' mollifiers $J_\epsilon\in OPS^{-\infty}$ and passing to the
limit as $\epsilon\to 0$ (cf. \cite{mis,tay}).

 \textbf{Step 1}. Estimates for the first component
$u$.

For $s\geq 2$, we calculate that (using \eqref{HS2} and taking
$h=0$)
\begin{eqnarray*}
\frac{d}{dt} \|   u \|_{H^s}^2 &=&
2 \langle \Lambda^s u_t, \Lambda^s   u \rangle\\
&=& -2 \langle \Lambda^s   ( u u_x ), \Lambda^s   u \rangle
+(\alpha+2) \langle \Lambda^s
\partial_x^{-1}   (u_x^2), \Lambda^s   u\rangle \\
&& +\kappa  \langle \Lambda^s \partial_x^{-1}   (\rho^2), \Lambda^s   u\rangle + 2 \langle \Lambda^s \partial_x^{-1} a, \Lambda^s   u\rangle\\
&:=& I_1+I_2+I_3+I_4.
\end{eqnarray*}
The first term can be estimated as in \cite{tay} by using the Kato-Ponce
estimate:
\begin{eqnarray*}
 |I_1| & =&2|\langle [\Lambda^s,  u] u_x, \Lambda^s   u\rangle
 + \langle   u  \Lambda^s u_x, \Lambda^s   u\rangle|
   \\
   &\leq & 2|\langle [\Lambda^s,  u] u_x, \Lambda^s   u\rangle|
  +|\langle
   u_x \Lambda^s   u,\Lambda^s   u\rangle|\\
   &\lesssim& \|[\Lambda^s,  u] u_x\|\|\Lambda^s   u\|+
   \| u_x\|_{L^\infty} \|\Lambda^s   u\|^2\\
   &\lesssim& \| u_x\|_{L^\infty} \|  u\|^2_{H^s}.
\end{eqnarray*}
For the second term, we make use of the Kato-Ponce commutator estimate (Lemma \ref{KP}):
\begin{eqnarray*}
 |I_2| &\leq & |\alpha+2|\|(u_x)^2\|_{H^{s-1}}\|\|  u\|_{H^s}\\
 &\lesssim& |\alpha+2|\|u_x\|_{L^\infty}\|u_x\|_{H^{s-1}}\|  u\|_{H^s}\\
 &\lesssim& |\alpha+2| \|u_x\|_{L^\infty}\|  u\|_{H^s}^2.
\end{eqnarray*}
Similarly, we can bound the third term involving the density $\rho$
by
\begin{eqnarray*}
|I_3| &\lesssim& |\kappa| \| \rho^2\|_{H^{s-1}}\|\| u\|_{H^s}
 \lesssim |\kappa| \|   \rho \|_{L^\infty}\|  \rho\|_{H^{s-1}}\|\|  u\|_{H^s}\\
&\lesssim& |\kappa| \|   \rho \|_{L^\infty}\| \rho\|^2_{H^{s-1}}+
|\kappa| \|   \rho \|_{L^\infty} \|  u\|^2_{H^s}.
\end{eqnarray*}

Since the quantity $a(t)$ only may depend on the time variable (cf.
\eqref{at} and Lemma \ref{conva}), we have $\Lambda^{2r}(ax)=ax$ for any
$r\geq 0$ being an integer. Therefore, we can bound $I_4$ as
follows:
\begin{eqnarray*}
 |I_4| &=& 2\left|\left\langle \Lambda^{s-2\left[\frac{s}{2}\right]}(ax), \Lambda^s   u\right\rangle\right| \lesssim
 \|ax\|_{H^1}\|  u \|_{H^s} \lesssim  |a| \|  u \|_{H^s}.
\end{eqnarray*}
A combination of the estimates for $I_1,...,I_4$ implies the
following inequality:
 \begin{eqnarray}
\frac{d}{d t} \|    u \|_{H^s}^2 &\lesssim & [(|\alpha+2|+1)\|
u_x \|_{L^\infty} +|\kappa| \|   \rho \|_{L^\infty}]
\|   u \|_{H^s}^2+ |a| \|  u \|_{H^s}\nonumber\\
&& + |\kappa| \|   \rho \|_{L^\infty}\| \rho\|^2_{H^{s-1}}.\label{u}
 \end{eqnarray}

 \textbf{Step 2}. Estimates for the second component $\rho$.

 We calculate that
 \begin{eqnarray}
 \frac{d}{dt}\|   \rho\|^2_{H^{s-1}}
 &=& 2\alpha \langle \Lambda^{s-1}  (u_x \rho), \Lambda^{s-1}   \rho\rangle-2 \langle \Lambda^{s-1}
 (\rho_x u ), \Lambda^{s-1}   \rho\rangle\nonumber\\
 &:=& J_1+J_2.\label{RHO}
 \end{eqnarray}
 The first term $J_1$ can be estimated like $I_2$ by Lemma
 \ref{AL}
 \begin{eqnarray}
  |J_1| &\lesssim& |\alpha| \|   u_x  \rho  \|_{H^{s-1}}\|  \rho\|_{H^{s-1}}
 \nonumber\\
 &\lesssim& |\alpha| (\|  u_x\|_{L^\infty}\|  \rho\|_{H^{s-1}}
 +\|  \rho \|_{L^\infty}\|  u_x \|_{H^{s-1}})\|  \rho\|_{H^{s-1}}\nonumber\\
 &\lesssim& |\alpha| (\|  u_x\|_{L^\infty}+\|  \rho
 \|_{L^\infty})(\|   u \|^2_{H^{s}}+\|  \rho\|^2_{H^{s-1}}).\label{J1}
 \end{eqnarray}
 Then we apply by the Kato-Ponce estimate (Lemma \ref{KP}) to
 $J_2$:
 \begin{eqnarray}
 |J_2|
 &=& 2|\langle [\Lambda^{s-1},  u] \rho_x, \Lambda^{s-1}   \rho\rangle
 + \langle   u  \Lambda^{s-1}   \rho_x, \Lambda^{s-1}   \rho\rangle|\nonumber\\
 &\leq & 2|\langle [\Lambda^{s-1},  u] \rho_x, \Lambda^{s-1}   \rho\rangle|+|\langle u_x  \Lambda^{s-1}    \rho, \Lambda^{s-1}   \rho\rangle|\nonumber\\
 &\lesssim& \|[\Lambda^{s-1},  u]  \rho_x\|\|\Lambda^{s-1}   \rho\|+\|u_x \|_{L^\infty}\|\Lambda^{s-1}    \rho\|^2\nonumber\\
 &\lesssim& \|u_x\|_{L^\infty}\|\Lambda^{s-2} \rho_x\|\|\Lambda^{s-1}   \rho\| +  \|\rho_x \|_{L^\infty} \|\Lambda^{s-1}   u \|\|\Lambda^{s-1}   \rho\|\nonumber\\
 && +\|u_x \|_{L^\infty}\|\Lambda^{s-1}
 \rho\|^2\nonumber\\
 &\lesssim& \|u_x \|_{L^\infty}\|\Lambda^{s-1}
 \rho\|^2+ \|\rho_x \|_{L^\infty} \|\Lambda^{s-1}   u \|\|\Lambda^{s-1}
 \rho\|.
 \label{rho1}
 \end{eqnarray}
 It follows from \eqref{RHO} -- \eqref{rho1} and the
H\"{o}lder inequality that
 \begin{eqnarray}
  \frac{d}{dt}\|   \rho\|^2_{H^{s-1}}
  \lesssim [(1+|\alpha|)\|u_x \|_{L^\infty}+ |\alpha|\|  \rho\|_{L^\infty}+ \|\rho_x \|_{L^\infty}]
   (\|  u \|_{H^s}^2+\|   \rho\|_{H^{s-1}}^2). \label{rhoa}
 \end{eqnarray}
 Combining \eqref{u} and \eqref{rhoa}, we can see that
 \begin{eqnarray}
 &&
  \frac{d}{d t} (\|    u \|_{H^s}^2+\|
  \rho\|^2_{H^{s-1}})\nonumber\\
 & \lesssim& \left[(1+|\alpha|+|\alpha+2|)\|u_x \|_{L^\infty}+
 (|\alpha|+|\kappa|)\|\rho\|_{L^\infty}+ \|\rho_x \|_{L^\infty}
 +|a|\right]\nonumber\\
&& \times  (\|u \|^2_{H^s}+\| \rho\|_{H^{s-1}}^2)+|a|.
 \end{eqnarray}
 From the definition of $a(t)$, we infer that
 \begin{equation}
 |a(t)|\leq
 \frac12\left|\kappa\right|\|\rho(t,\cdot)\|_{L^\infty}^2+\frac12
 \left|\alpha+2\right|\| u_x(t,\cdot)\|_{L^\infty}^2.\label{BDa1}
 \end{equation}
 \begin{rem}
 We note that when $\alpha=-1$, $|a(t)|$ is a constant that can be controlled by $\|u_x^0\|,  \|\rho^0\|$ and $|\kappa|$.
 \end{rem}
 Under the assumption \eqref{Asm}, for $t\in [0,T)$, it
 holds
 \begin{eqnarray}
 && \frac{d}{d t} (\|u \|_{H^s}^2+\|  \rho\|^2_{H^{s-1}})\nonumber\\
 &\lesssim&
 \left[(1+|\alpha|+|\alpha+2|+|\kappa|) M+
 (|\alpha+2|+|\kappa|)M^2\right]
  (\| u \|^2_{H^s}+\| \rho\|_{H^{s-1}}^2) \nonumber\\
  && +(|\alpha+2|+|\kappa|)M^2. \nonumber
 \end{eqnarray}
 By the Gronwall inequality, we see that $\|(u,\rho)^{tr}\|_{H^s\times H^{s-1}}$ is bounded
 for $t\in [0,T)$. The proof is complete.
\end{proof}

\section{Blow-up scenarios for $(\alpha,\kappa)\in \{-1\}\times \mathbb{R}_+$}
 \setcounter{equation}{0}

 In Section \ref{Per}, we have shown a persistence result for all
 $(\alpha, \kappa)\in \mathbb{R}\times\mathbb{R}$. Concerning the interested case of $(\alpha,\kappa)\in \{-1\}\times
 \mathbb{R}_+$, we can prove the precise blow-up scenarios for regular
 solutions.

\begin{lem}\label{H2n}
Suppose that $(\alpha,\kappa)\in \{-1\}\times \mathbb{R}_+$. Given
any $z^0=(u^0, \rho^0)^{tr}\in H^s(\mathbb{S})\times
H^{s-1}(\mathbb{S})$, $s\geq 2$. For the solution $z=(u, \rho)^{tr}$
of system \eqref{hssys} corresponding to $z^0$, we have
 \begin{equation}
 \|\rho(t)\|^2+\|u_x(t)\|^2\leq C(\|\rho^0\|, \|u^0_x\|),\quad \forall \
 t>0.\label{ux1}
 \end{equation}
 Moreover,
\begin{equation}
\|u(t)\|\leq C(T,\|u^0\|_{H^1},\|\rho^0\|), \quad t\in
[0,T).\label{uL2}
\end{equation}
\end{lem}
\begin{proof}
Estimate \eqref{ux1} follows from Lemma \ref{conva} and the facts that $\alpha=-1$,
$\kappa>0$. For the proof of \eqref{uL2}, we refer to \cite[pp.
653]{Wun09}.
\end{proof}

Let $u(t,x)$ be the solution of \eqref{hssys}. We consider the initial value problem for the Lagrangian flow
map:
 \begin{equation} \partial_t \; \varphi(t,x)=u(t,\varphi(t,x)), \quad
\varphi(0,x)=x.\label{phi}
 \end{equation}
 It is well-known that  (cf. e.g., \cite{mis}) the following lemma is valid.
 \begin{lem}
Let $u\in C([0,T);H^s)\cap C^1([0,T);H^{s-1}), s\geq 2$. Then
problem \eqref{phi} admits a unique solution $x\in
C^1([0,T)\times\mathbb{S};\mathbb{S})$. Moreover,
$\{\varphi(t,\cdot)\}_{t\in [0,T)}$ is a family of
orientation-preserving diffeomorphisms on the circle $\mathbb{S}$
and
\begin{equation}
 \varphi_x(t,x)=e^{\int_0^t u_x(s, \varphi(s,x)) ds}>0,\quad (t,x)\in
[0,T)\times \mathbb{S}.\label{LFM}
\end{equation}
 \end{lem}
Moreover, if $\alpha=-1$, in analogy to \cite[Lemma 3.4]{ELY07}, we can show
that
\begin{lem}
Suppose that $(\alpha,\kappa)\in \{-1\}\times \mathbb{R}_+$. Given
any $z^0=(u^0, \rho^0)^{tr}\in H^2(\mathbb{S})\times
H^1(\mathbb{S})$. Let $z=(u, \rho)^{tr}$ be the solution to system
\eqref{hssys} corresponding to $z^0$ on $[0,T)$. We have
\begin{equation}
\rho(t,\varphi(t,x))\varphi_x(t,x)=\rho^0(x), \quad \forall (t,x)\in
[0,T)\times \mathbb{S}.\label{LFMR}
\end{equation}
Moreover, if there exists $M_1 > 0$ such that $u_x(t, x)\geq- M_1$
for all $(t, x) \in [0, T )\times\mathbb{S}$, then
\begin{equation}
\|\rho(t,\cdot)\|_{L^\infty}=\|\rho(t,\varphi(t,\cdot)\|_{L^\infty}\leq
e^{M_1T}\|\rho^0(\cdot)\|_{L^\infty},\quad t\in [0,T).\label{rhoinf}
\end{equation}
\end{lem}

\begin{thm}\label{BB1}
Suppose that $(\alpha, \kappa)\in \{-1\}\times\mathbb{R}_+$. For any
$z^0=(u^0, \rho^0)^{tr}\in H^2(\mathbb{S})\times H^{1}(\mathbb{S})$,
let $T$ be the maximal existence time of the solution $z=(u,
\rho)^{tr}$ to \eqref{hssys} corresponding to the initial datum $z^0$. Then the
solution blows up in finite time if and only if
\begin{equation}
 \liminf_{t\to T^-}\{\inf_{x\in\mathbb{S}}
 u_x(t,x)\}=-\infty.\label{blow1}
\end{equation}
\end{thm}
\begin{proof}
Consider the equation describing the dynamics of $\rho(t, x)$ in
\eqref{hssys}, and differentiate it once in space
 \begin{equation}
\rho_{xt} + u\rho_{xx}+2u_x\rho_{x}+u_{xx}\rho=0.\label{drho}
 \end{equation}
Multiplying the first equation in \eqref{hssys} by $u_{xx}$, and
\eqref{drho} by $\rho_x$, upon adding the resultants together, we
deduce that
\begin{equation}
\frac12\partial_t(u_{xx}^2+\rho_x^2)+2u_x(u_{xx}^2+\rho_x^2)+\frac12
u\partial_x(u_{xx}^2+\rho_x^2)+(1-\kappa)\rho\rho_xu_{xx}=0.\label{high}
\end{equation}
Integrating in space and using the periodic boundary conditions,
we have
 \begin{eqnarray}
&&
\frac{d}{dt}(\|u_{xx}\|^2+\|\rho_x\|^2)\nonumber\\
&=&-3\int_\mathbb{S} u_x (u_{xx}^2+\rho_x^2) dx+
(1-\kappa)\int_\mathbb{S}\rho\rho_xu_{xx}
dx\nonumber\\
&\leq&-3\int_\mathbb{S} u_x (u_{xx}^2+\rho_x^2)
dx+|1-\kappa|\|\rho\|_{L^\infty}\int_\mathbb{S} (u_{xx}^2+\rho_x^2)
dx.\label{h2r1}
\end{eqnarray}
Assume that there exists $M_1 > 0$ such that
\begin{equation}
u_x(t,x)\geq -M_1, \quad \forall \ (t,x)\in
[0,T)\times\mathbb{S}.\label{bounduxa}
\end{equation}
Then it follows from \eqref{rhoinf} and \eqref{h2r1} that
\begin{eqnarray}
 \frac{d}{dt}\left(\|u_{xx}\|^2+\|\rho_x\|^2\right)\leq
 \left(3M_1+|1-\kappa|e^{M_1T}\|\rho^0\|_{L^\infty}\right)\left(\|u_{xx}\|^2+\|\rho_x\|^2\right).
\end{eqnarray}
By Gronwall's inequality we have
 \begin{equation}
 \|u_{xx}(t)\|^2+\|\rho_x(t)\|^2\leq \left(\|u^0_{xx}\|^2+\|\rho^0_x\|^2\right)
 e^{\left(3M_1+|1-\kappa|e^{M_1T}\|\rho^0\|_{L^\infty}\right)t},\quad
 t\in [0,T).
 \end{equation}
This and Lemma \ref{H2n} ensure that the solution $z$ does not blow
up in finite time.

On the other hand, by Theorem \ref{T1}, we see that if \eqref{blow1}
holds, then the solution will blow up in finite time. The proof is
complete.
\end{proof}

\begin{lem}
\label{hb1} Suppose that $(\alpha, \kappa)\in \{-1\}\times\mathbb{R}_+$.
Let $z^0=(u^0, \rho^0)^{tr}\in H^s(\mathbb{S})\times
H^{s-1}(\mathbb{S})$, $s>2$, and let $T$ be the maximal existence
time of the solution $z=(u, \rho)^{tr}$ to \eqref{hssys} with the
initial datum $z^0$. If there exist two constants $M_1, M_2>0$ such that
\begin{equation}
u_x(t,x)\geq -M_1, \quad \|\rho_x(t,\cdot)\|_{L^\infty}\leq
M_2,\quad \forall \ (t,x)\in [0,T)\times\mathbb{S},\label{blow3}
\end{equation}
then $\|z(t, \cdot)\|_{H^s\times H^{s-1}}$ will not blow up in finite time.
\end{lem}
\begin{proof}
Under our current assumption \eqref{blow3}, it follows from the argument in
Theorem \ref{BB1} that $\|z\|_{H^2\times H^{1}}$ is bounded for all
$t\in [0,T)$. By Sobolev's embedding theorem, we can see that
$\|u_x\|_{L^\infty}$ and $\|\rho\|_{L^\infty}$ are also bounded.
Then our conclusion  easily follows from Theorem \ref{T1}.
\end{proof}

Now we show a first precise blow-up scenario for sufficiently
regular solutions:

\begin{thm}\label{B1}
Suppose that $(\alpha, \kappa)\in \{-1\}\times\mathbb{R}_+$. For any
$z^0=(u^0, \rho^0)^{tr}\in H^s(\mathbb{S})\times
H^{s-1}(\mathbb{S})$, $s>\frac52$, let $T$ be the maximal
existence time of the solution $z=(u, \rho)^{tr}$ to \eqref{hssys}
with initial datum $z^0$. Then the corresponding solution blows
up in finite time if and only if
\begin{equation}
 \liminf_{t\to T^-}\{\inf_{x\in\mathbb{S}}
 u_x(t,x)\}=-\infty,\quad \text{or} \quad \limsup_{t\to T^-}\|\rho_x(t,\cdot)\|_{L^\infty}=+\infty.\label{blow2}
\end{equation}
\end{thm}
\begin{proof}
By Lemma \ref{hb1}, we can see that if there exist $M_1, M_2>0$ such
that assumption \eqref{blow3} are satisfied, then $\|z\|_{H^s\times H^{s-1}}$
will not blow up in finite time. On the other hand, by Theorem
\ref{T1} and Sobolev's embedding theorem, we see that if
\eqref{blow2} holds, then the solution will blow up in finite time. The proof is complete.
\end{proof}

We note that a blow-up scenario similar to Theorem \ref{B1} that
involves the condition \eqref{blow2} was obtained for regular
solutions to a two-component Camassa-Holm equations (cf. e.g.,
\cite{ELY07}). Later in \cite{ZL09}, the authors obtained an
improved blow-up scenario that only needs the condition on one of
the component (i.e., \eqref{blow1} for $u$). In what follows, we
derive an improved blow-up scenario for our two-component
Hunter-Saxton system, which shows that the assumption \eqref{blow1},
is actually enough to determine wave breaking of the regular
solutions ($s>\frac52$) in finite time. The key observation is that
the quantity $\|\rho_x\|_{L^\infty}$ can be controlled by the
lower-bound of $u_x$.

\begin{thm}\label{B2}
Suppose that $(\alpha, \kappa)\in \{-1\}\times\mathbb{R}_+$. For any
$z^0=(u^0, \rho^0)^{tr}\in H^s(\mathbb{S})\times
H^{s-1}(\mathbb{S})$, $s>\frac52$, let $T$ be the maximal existence
time of the solution $z=(u, \rho)^{tr}$ to \eqref{hssys} with
initial datum $z^0$. Then the corresponding solution blows up in
finite time if and only if \eqref{blow1} holds.
\end{thm}
\begin{proof}
Using the Lagrangian flow map, we set
$$ M(t,x)=u_x(t, \varphi(t,x)),\quad \gamma(t,x)=\rho(t,
\varphi(t,x)), $$
 $$ N(t,x)=u_{xx}(t, \varphi(t,x)),\quad \varpi(t,x)=\rho_x(t,\varphi(t,x)). $$
It follows from \eqref{high} that
 \begin{equation}
\partial_t(N^2+\varpi^2)+4M(N^2+\varpi^2)+2(1-\kappa)\gamma\varpi N=0.\label{high1}
\end{equation}
Then we have
 \begin{equation}
\partial_t(N^2+\varpi^2)\leq
(-4M+|1-\kappa|\|\gamma\|_{L^\infty})(N^2+\varpi^2). \label{highli}
 \end{equation}
 Assume that there exists $M_1 > 0$ such that \eqref{bounduxa} holds.
Then it follows from \eqref{rhoinf} and \eqref{highli} that
\begin{eqnarray}
 \frac{d}{dt}\left(N^2+\varpi^2\right)\leq
 \left(4M_1+|1-\kappa|e^{M_1T}\|\rho^0\|_{L^\infty}\right)\left(N^2+\varpi^2\right).
\end{eqnarray}
By Gronwall's inequality and the Sobolev embedding theorem
($s>\frac52$), for all $(t,x)\in [0,T)\times \mathbb{S}$, we have
 \begin{eqnarray}
 && N(t,x)^2+\varpi(t,x)^2\nonumber\\
 &\leq& \left(\|u^0_{xx}\|_{L^\infty}^2+\|\rho^0_x\|_{L^\infty}^2\right)
 e^{\left(4M_1+|1-\kappa|e^{M_1T}\|\rho^0\|_{L^\infty}\right)t}\nonumber\\
 &\leq& C\left(\|u^0\|_{H^s}^2+\|\rho^0\|_{H^{s-1}}^2\right)
 e^{\left(4M_1+|1-\kappa|e^{M_1T}\|\rho^0\|_{L^\infty}\right)t}.
 \end{eqnarray}
In particular, this implies that there exists a constant $M_2>0$
such that
\begin{equation}
\|\rho_x(t,\cdot)\|_{L^\infty}\leq M_2,\quad \forall \ t\in [0,T).
\end{equation}
Thus, the condition \eqref{blow3} made in Lemma \ref{hb1} is now
satisfied, and as a result, we conclude that $\|z(t,
\cdot)\|_{H^s\times H^{s-1}}$ will not blow up in finite time.

On the other hand, by Theorem \ref{T1}, we see that if \eqref{blow1}
holds, then the solution will blow up in finite time. The proof is
complete.
\end{proof}

\section{Global existence for $(\alpha,\kappa)\in \{-1\}\times \mathbb{R}_+$}
\setcounter{equation}{0}

We remark here that the assumption $\eqref{Asm}$ is somewhat artificial:
There is no reason to assume that
$\|\rho_x(t, \cdot)\|_{L^\infty}$ actually stays bounded in time
(but note that this assumption was also made in \cite{ELY07} for the 2-component Camassa-Holm equations).
It turns out, however, that we can dispense with \eqref{Asm} if we impose some sign condition on the initial datum $\rho^0$. Our results show that if $\rho^0(x)$ keeps its sign for all $x\in \mathbb{S}$, then existence of global solutions to system \eqref{hssys} will be guaranteed for $(\alpha,\kappa)\in \{-1\}\times \mathbb{R}_+$. Besides, in the previous work \cite{Wun10}, a smallness condition on the quantity $\|u_x^0\|^2+\kappa\|\rho^0\|^2$ was required to obtain the (global-in-time) lower-order estimate of the solutions. In what follows, we improve the former results by showing that only the sign condition of the initial data can ensure the existence of regular solutions of our system.

\subsection{Global existence in $H^2\times H^1$}

\begin{thm} \label{-1H2}
Suppose that $(\alpha,\kappa)\in \{-1\}\times \mathbb{R}_+$. Given
any $z^0=(u^0, \rho^0)^{tr}\in H^2(\mathbb{S})\times
H^1(\mathbb{S})$, if we further assume that
 \begin{equation} \rho^0(x)>0,\quad \text{or} \quad \rho^0(x)<0,  \label{condi} \end{equation}
then the solution $z=(u, \rho)^{tr}$ to system \eqref{hssys}
corresponding to $z^0$ is global.
\end{thm}
\begin{proof}
It suffices to get some uniform \emph{a priori} estimates for the
solution $(u, \rho)^{tr}$.

 \textbf{Step 1}. Estimates for $\|\rho\|_{L^\infty}$ and
$\|u_x\|_{L^\infty}$.

Using the Lagrangian flow map, we set
 $$ M(t,x)=u_x(t, \varphi(t,x)),\quad \gamma(t,x)=\rho(t,
\varphi(t,x)).$$
 Then we have (cf. \cite{ci,Wun09,Wun10})
 \begin{equation}
 M_t(t,x)=-\frac12 M^2(t,x)+\frac\kappa2 \gamma(t,x)^2+a,\quad
 \gamma_t(t,x)=-M(t,x)\gamma(t,x).
 \end{equation}
By the assumptions \eqref{condi}, we infer from \eqref{LFM} and
\eqref{LFMR} that if $\gamma(0,x)>0$ (or $\gamma(0,x)<0$) then
$\gamma(t,x)>0$ (or  $\gamma(t,x)<0$) for $t\in [0,T)$. Thus, we can
construct the following strictly positive auxiliary function (cf.
\cite{ci,Wun09})
 $$ w(t,x):=
\kappa\gamma(0,x)\gamma(t,x)+\frac{\gamma(0,x)}{\gamma(t,x)}[1+M(t,x)^2].$$
 Computing the evolution of $w$, we get
\begin{eqnarray*}
\partial_t \;w(t,x) &=& \kappa\gamma(0,x) \; \partial_t \gamma(t,x) - \frac{\gamma(0,x)}{\gamma(t,x)^2} \; \partial_t
\gamma(t,x)\left[1 + M(t,x)^2\right] \\
&&+ 2 \frac{\gamma(0,x)}{\gamma(t,x)}M(t,x)\;\partial_t M(t,x)\\
&=& -\kappa\gamma(0,x)\gamma(t,x)M(t,x) + \frac{\gamma(0,x)}{\gamma(t,x)^2}\gamma(t,x)M(t,x)\left(1 + M^2(t,x)\right) \\
&&- \frac{\gamma(0,x)}{\gamma(t,x)}M(t,x)^3 +\kappa
\frac{\gamma(0,x)}{\gamma(t,x)}M(t,x)\gamma(t,x)^2
+ 2a\frac{\gamma(0,x)}{\gamma(t,x)}M(t,x)\\
&=& (1+ 2a)\frac{\gamma(0,x)}{\gamma(t,x)}M(t,x).
\end{eqnarray*}
 The last quantity can be estimated by
 \begin{eqnarray*}
(1 + 2 |a|) \frac{\gamma(0,x)}{\gamma(t,x)} |M(t,x)|
&\le & (1 + 2|a|) \frac{\gamma(0,x)}{\gamma(t,x)} (1 + M^2(t,x))\\
&\le& (1 + 2|a|)w(t,x).
 \end{eqnarray*}
By Gronwall's inequality, we obtain
\begin{equation}
 \label{boundux}
w(t,x) \leq w(0,x)\;e^{(1+2|a|)t}, \quad t\in [0,T),
\end{equation}
 which together with \eqref{conaa} implies the following estimate
\begin{eqnarray}
\|\rho(t)\|_{L^\infty}+\|u_x(t)\|_{L^\infty}&\leq&
C(T,\|\rho^0\|_{L^\infty},\|u^0_x\|_{L^\infty}
)\nonumber\\
&\leq& C(T,\|\rho^0\|_{H^1},\|u^0\|_{H^2} ).\label{inf}
\end{eqnarray}

\textbf{Step 2}. Estimates for $\|u\|_{H^2}$ and $\|\rho\|_{H^1}$.

It follows from \eqref{h2r1} that
\begin{eqnarray}
 \frac{d}{dt}\left(\|u_{xx}\|^2+\|\rho_x\|^2\right)\leq
 \left(3\|u_x\|_{L^\infty}+|1-\kappa|\|\rho\|_{L^\infty}\right)\left(\|u_{xx}\|^2+\|\rho_x\|^2\right).
\end{eqnarray}
By Gronwall's inequality we have
$$ \|u_{xx}(t)\|^2+\|\rho_x(t)\|^2\leq e^{(3\|u_x\|_{L^\infty}+|1-\kappa|\|\rho\|_{L^\infty})t}(\|u^0_{xx}\|^2+\|\rho^0_x\|^2),$$
 which together with \eqref{inf} and Lemma \ref{H2n} yields that
\begin{equation}
\|u(t)\|_{H^2}^2+\|\rho(t)\|_{H^1}^2\leq
C(T,\kappa, \|\rho^0\|_{H^1},\|u^0\|_{H^2}),\quad t\in [0,T).\label{H2}
\end{equation}
\end{proof}

\begin{rem}
We notice that, for $\kappa=1$, one only needs a bound on $\|u_x(t,.)\|_{L^\infty}$
to get the uniform estimate \eqref{H2}. Besides, in contrast with \cite[Proposition 6.1]{Wun09}, we have shown that
in order to have global existence in $H^2\times H^1$, one does not
need to impose certain smallness assumptions on the initial data, and it only requires that $\rho^0$ is strictly nonzero (cf. \eqref{condi}). This follows from an idea of \cite{guo}.
\end{rem}

\subsection{Global existence in $H^s\times H^{s-1}$ }

\begin{thm}\label{T2}
Suppose that $(\alpha,\kappa)\in \{-1\}\times \mathbb{R}_+$. Given
any $z^0=(u^0, \rho^0)^{tr}\in H^s(\mathbb{S})\times
H^{s-1}(\mathbb{S})$, $s>
 \frac52$ that satisfies \eqref{condi}. Assume that $T$ is the existence time of the solution $z=(u,
 \rho)^{tr}$ to system \eqref{hssys} corresponding to $z^0$. Then $\|z(t, \cdot)\|_{H^s\times H^{s-1}}$ is bounded on $[0,T)$.
\end{thm}
\begin{proof}
 The key point is how to estimate $
 \|\rho_x \|_{L^\infty}$, which is the main difficulty for the second term in the last line of
 \eqref{rho1}.

 By the Gagliardo-Nirenberg inequality, we see that
  \begin{equation}
  \|  \rho_x\|_{L^\infty}\lesssim\|   \rho\|_{H^{s-1}}^p\|   \rho\|_{L^\infty}^{1-p}+\|\rho\|_{L^\infty},\quad
  p  =\frac{1}{s-\frac32}\nonumber
  \end{equation}
 and
 \begin{eqnarray}
 \|\Lambda^{s-1}   u \|&\lesssim& \|\Lambda^{s-2}  u_x \|+ \|  u \|_{H^1}\nonumber\\
 &\lesssim& \|\Lambda^{s-1} u_x \|^q\|u_x \|^{1-q}_{L^\infty}+\|u_x \|_{L^\infty}
 + \|  u \|_{H^1},\ \ q=\frac{s-\frac52}{s-\frac32}.\nonumber
 \end{eqnarray}
 \begin{rem} We notice that $p+q=1$ and $p>0$ for all $s>\frac32$.
 If we take $s> \frac52$, then $p\in (0,1)$ and as a result, $q\in
 (0,1)$. We note that $s=\frac52$ is an exceptional case of the Gagliardo-Nirenberg
 inequality.
 \end{rem}
 In what follows, we assume that $s> \frac52$. The second term in the last line of
 \eqref{rho1} can be estimated as follows
 \begin{eqnarray}
  && \|\rho_x \|_{L^\infty} \|\Lambda^{s-1}   u \|\|\Lambda^{s-1}   \rho\|\nonumber\\
  &\lesssim& (\|   \rho\|_{L^\infty}^{1-p} \|   \rho\|_{H^{s-1}}^{1+p}+\|\rho\|_{L^\infty}\|\rho\|_{H^{s-1}})
  (\|  u \|_{H^s}^{1-p}\| u_x \|^{p}_{L^\infty}+ \|u_x \|_{L^\infty}+\|  u \|_{H^1})\nonumber\\
  &\lesssim& \|  \rho\|_{L^\infty}^{1-p}\|u_x \|^{p}_{L^\infty} (\|  u \|_{H^s}^2+\|   \rho\|_{H^{s-1}}^2)+\|   \rho\|_{L^\infty}^{1-p} \|   \rho\|_{H^{s-1}}^{1+p}(\|u_x \|_{L^\infty}+\|  u
  \|_{H^1})\nonumber\\
  && +\|\rho\|_{L^\infty}\|\rho\|_{H^{s-1}}(\|u\|_{H^s}+\|u_x\|_{L^\infty}+\|u\|_{H^1})\nonumber\\
  &\lesssim& (\|  \rho\|_{L^\infty} + \|u_x \|_{L^\infty}) (\|  u \|_{H^s}^2+\|   \rho\|_{H^{s-1}}^2)\nonumber\\
  && +(\|   \rho\|_{L^\infty}^2+ \|   \rho\|_{H^{s-1}}^2)(\|u_x \|_{L^\infty}+\|  u \|_{H^1}).\label{d}
 \end{eqnarray}
 As a result, we conclude from \eqref{rho1} and \eqref{d} a new estimate for $J_2$:
 \begin{eqnarray}
 |J_2|  &\lesssim& (\|  \rho\|_{L^\infty} + \|u_x \|_{L^\infty}) (\|  u \|_{H^s}^2+\|   \rho\|_{H^{s-1}}^2)\nonumber\\
  && +(\|u_x \|_{L^\infty}+\|  u \|_{H^1})(\|   \rho\|_{L^\infty}^2+ \|   \rho\|_{H^{s-1}}^2).\label{J2}
 \end{eqnarray}
 It follows from \eqref{RHO}, \eqref{J1}, \eqref{J2} and the Young inequality that
 \begin{eqnarray}
 \frac{d}{dt}\|   \rho\|^2_{H^{s-1}}
 &\lesssim& (\|u_x \|_{L^\infty}+\|  \rho\|_{L^\infty}) (\|  u \|_{H^s}^2+\|   \rho\|_{H^{s-1}}^2)\nonumber\\
 && +(\|u_x \|_{L^\infty}+\|  u \|_{H^1})(\|   \rho\|_{L^\infty}^2+ \|   \rho\|_{H^{s-1}}^2). \label{rho}
 \end{eqnarray}
Combining \eqref{u} and \eqref{rho}, we can see that
 \begin{eqnarray}
 \frac{d}{d t} (\|    u \|_{H^s}^2+\|   \rho\|^2_{H^{s-1}})
 & \lesssim& \left[ \|u_x \|_{L^\infty}+ (|\kappa|+1)\|\rho\|_{L^\infty}+\|  u \|_{H^1}+|a|\right]\nonumber\\
&& \times  (\|  u \|^2_{H^s}+\|  \rho\|_{H^{s-1}}^2)+(\|u_x \|_{L^\infty}+\|  u
\|_{H^1})\|  \rho\|_{L^\infty}^2+|a|,\nonumber
 \end{eqnarray}
 Since  $a$ is now a constant (cf. Lemma \ref{conva}), and $\|u(t)\|_{H^1}$ is bounded
 by $\|u^0\|_{H^1}, \|\rho^0\|$ and $T$ (cf. Lemma \ref{H2n}), it holds
\begin{eqnarray}
 \frac{d}{d t} (\|  u \|_{H^s}^2+\|  \rho\|^2_{H^{s-1}})
 & \leq& C_1 (\|  u_x\|_{L^\infty}+\| \rho\|_{L^\infty}+1)(\| u \|_{H^s}^2+\| \rho\|_{H^{s-1}}^2)\nonumber\\
 && +C_2\| \rho\|_{L^\infty}^2+C_3\|u_x \|_{L^\infty}\| \rho\|_{L^\infty}^2+C_4,\label{PER}
 \end{eqnarray}
 where $C_i (i=1,...,4)$ are positive constants that may depend on $\| u^0\|_{H^1}$, $\|\rho^0\|$ and $\kappa$, moreover,
 $C_1, C_2, C_3$ may also depend on $T$.

 We have known that if $(u^0, \rho^0)\in H^2\times H^1$, there exists a constant $K>0$ such
 that (cf. \eqref{inf})
 $$\| u_x(t)\|_{L^\infty}+\| \rho(t)\|_{L^\infty}\leq K,\quad \forall\ t\in [0,T).$$
 By the Gronwall inequality, we obtain that for all $t\in [0,
 T)$:
 \begin{equation}
 \|  u (t)\|_{H^s}^2+\|  \rho(t) \|^2_{H^{s-1}}\leq e^{C_1(K+1)t}(\|  u^0 \|_{H^s}^2
 +\|\rho^0 \|^2_{H^{s-1}})+\frac{C_2K^2+C_3K^3+C_4}{C_1(K+1)}.
 \end{equation}
 \end{proof}

 \begin{rem}
 Comparing Theorem \ref{-1H2} and Theorem \ref{T2}, we can see that there is a gap, namely, $s\in (2,\frac52]$.
One possible reason is that, for the case $s=2$, there is a special
structure in the derivation of
the evolution \eqref{high}.
However,
for general $s$, one has to deal with the term \eqref{d} in $J_2$ by
using proper interpolation inequalities, which excludes the case $s\in
(2,\frac52]$.
 \end{rem}

 \section{Global existence for $(\alpha,\kappa)\in \{0\}\times \mathbb{R}_+$}
\setcounter{equation}{0}

From Lemma \ref{conva}, we see that a very important property for
the case $(\alpha,\kappa)\in \{-1\}\times \mathbb{R}_+$ is the
conservation law for the quantity $a(t)$, which can be bounded by $\|u^0_x\|$,
$\|\rho^0\|$. However, this nice property may be lost for other choices
of $\alpha$ (cf. Lemma \ref{conva}). This fact leads to a different procedure to prove that
solutions exist globally.

 \begin{thm} \label{T3}
Suppose that $(\alpha,\kappa)\in \{0\}\times \mathbb{R}_+$. Given
any $z^0=(u^0, \rho^0)^{tr}\in H^s(\mathbb{S})\times
H^{s-1}(\mathbb{S})$, $s\geq 3$, we assume that $\rho^0$ satisfies
the sign condition \eqref{condi}, and $T$ is the existence time of
the solution $z=(u,
 \rho)^{tr}$ to system \eqref{hssys} corresponding to $z^0$.
 Then $\|z(t, \cdot)\|_{H^s\times H^{s-1}}$ is bounded on $[0,T)$.
\end{thm}
\begin{proof}
\textbf{Step 1}. Estimates for $\|\rho\|_{L^\infty}$,
$\|u_x\|_{L^\infty}$.

Now for $\alpha=0$, we no longer have the conservation law for $a(t)$
(cf. Lemma \ref{conva}). As a result, we lose the control of
$\|\rho\|, \|u_x\|$, in contrast with the case $\alpha=-1$ (see Lemma
\ref{H2n}). Fortunately, however, the equation for $\rho$ now is just a
transport equation, which implies that
 \begin{equation}
  \|\rho(t,\cdot)\|_{L^\infty}=\|\rho^0\|_{L^\infty}\leq
 C(\|\rho^0\|_{H^2}), \quad t\in [0, T).\label{rhobd1}
 \end{equation}
 Besides, it follows from  \cite[Proposition 4.1]{Wun10} that for $t\in
 [0,T)$,
 \begin{eqnarray}
  \sup_{x\in \mathbb{S}}u_x(t,x)&\leq& C(\|u^0\|_{H^3}, \|\rho^0\|_{H^2},\kappa,
 T),\label{supux} \\
  \|u_x(t)\|&\leq& C(\|u^0\|_{H^3}, \|\rho^0\|_{H^2},\kappa,
 T).\label{Ux1}
 \end{eqnarray}
 Multiplying \eqref{HS2} by $u$, integrating over the circle, we
 have
 \begin{eqnarray}
 && \frac12 \frac{d}{dt} \|u\|^2\nonumber\\
  &=& \int_\mathbb{S}
 u\partial_x^{-1}\left(\frac12 \rho^2+u_x^2+a(t)\right) dx +h(t)\int_\mathbb{S}
 u dx\nonumber\\
 &\leq&\frac12|h(t)|+ \frac12(1+|h(t)|)\int_\mathbb{S} u^2dx +
 \frac12 \left[\int_\mathbb{S} \left(\frac12 \rho^2+u_x^2+|a(t)|\right)
 dx\right]^2.
 \end{eqnarray}
 It follows from \eqref{rhobd1}, \eqref{Ux1} that
 \begin{equation} |a(t)|\leq  \frac{\kappa}{2}\|\rho(t)\|_{L^\infty}^2+ \|u_x(t)\|^2\leq C(\|u^0\|_{H^3},
 \|\rho^0\|_{H^2},\kappa,
 T), \quad \forall t\in [0, T),\label{abd}
 \end{equation}
 which yields
 \begin{equation}
  \frac{d}{dt} \|u\|^2\leq C(1+\|u\|^2)+C',
 \end{equation}
 where $C$ is a constant depending on $h(t)$ and $C'$ is a
 constant depending on $ \|u^0\|_{H^3}$, $\|\rho^0\|_{H^2}$, $T$ and $\kappa$.
 By the Gronwall inequality and \eqref{Ux1}, we see that
 \begin{equation}
 \|u(t)\|_{H^1}\leq C(\|u^0\|_{H^3}, \|\rho^0\|_{H^2},\kappa,
 T),\quad t\in [0,T). \label{uH1aa}
 \end{equation}

 Recalling the functions $M, \gamma$ introduced in the
 proof of Theorem \ref{-1H2}, in our present case we have
 \begin{equation}
 \partial_t M(t,x)=
 \frac\kappa 2\gamma(t,x)^2+a(t),\quad \partial_t\gamma (t,x)=0.
 \end{equation}
 Then we compute the time derivative of $$\tilde w(t,x)=\kappa \gamma (0,x)^2+ (1+M(t,x)^2)$$ such that
 \begin{eqnarray}
 \partial_t \tilde w&=&2M\partial_t M=
 (\kappa\gamma^2+2a)M\nonumber\\
 &\leq& \frac\kappa2\gamma^2(1+M^2)+|a|(1+M^2)\leq
 \left(\frac\kappa 2\gamma^2+|a|\right)\tilde w.
 \end{eqnarray}
It follows from the Gronwall inequality, \eqref{rhobd1},
\eqref{abd}, and the definition of $\tilde w$ that
 \begin{equation}
 \|u_x(t,\cdot)\|_{L^\infty}\leq C(\|u^0\|_{H^3}, \|\rho^0\|_{H^2},
 \kappa, T),\quad t\in [0,T).\label{Uxin1}
 \end{equation}

 \textbf{Step 2}. Estimates of $\|u\|_{H^s}$, $\|\rho\|_{H^{s-1}}$ ($s\geq 3$).

 Since now $s\geq 3>\frac52$, we can check the calculations step by step
 in the proof of Theorem \ref{T1} and Theorem \ref{T2} to see that the inequality \eqref{PER} still holds
 in the present situation (i.e., $(\alpha, \kappa)\in \{0\}\times \mathbb{R}_+$).
 However, due to \eqref{uH1aa}, the constants $C_1,...,C_4$ in
 \eqref{PER} depend on $\|u^0\|_{H^3}, \|\rho^0\|_{H^2}$, $\kappa$ and $ T$.
 By the Gronwall inequality and \eqref{rhobd1}, \eqref{Uxin1}, we obtain that for all $s\geq 3$ and $t\in [0,
 T)$:
 \begin{equation}
 \|u(t)\|_{H^s}^2+\|  \rho(t) \|^2_{H^{s-1}}\leq C(\|u^0\|_{H^s},
 \|\rho^0\|_{H^{s-1}}, \kappa, T).
 \end{equation}
 The proof is complete.
 \end{proof}
 \begin{rem}
 In Theorem \ref{T3}, we assumed that $s\geq 3$. This is because in
order to obtain the estimate \eqref{supux}, one has to make use
 of an abstract lemma due to {\sc Constantin \& Escher} \cite{ce} which requires that $u_x \in C^1([0, T];
 H^1)$, i.e., $z^0\in H^3\times H^2$ (cf. e.g., \cite{Wun10}).
 \end{rem}



\section*{Acknowledgements}
A part of this paper was written during the first author's visit to
RIMS in Kyoto University, whose hospitality is gratefully
acknowledged. The first author was supported by Natural Science
Foundation of Shanghai 10ZR1403800.
\par
The second author was supported
by JSPS Postdoctoral Fellowship P09024. He is grateful to Z. Guo for
fruitful discussions, and to Y. Zhou for his hospitality during the
Chinese-German Conference on Fluid and Gas Dynamics at Zhejiang
Normal University in May 2010.

\end{document}